\documentclass[12pt,a4paper]{amsart}
\usepackage{amsmath, amssymb, amsfonts, enumerate}
\usepackage[all]{xy}
\usepackage{amscd}

\newtheorem{theorem}{Theorem}[section]
\newtheorem{lemma}[theorem]{Lemma}
\newtheorem{corollary}[theorem]{Corollary}
\newtheorem{proposition}[theorem]{Proposition}

 \theoremstyle{definition}
 
 \newtheorem{remark}[theorem]{Remark}
\newtheorem{remarks}[theorem]{Remarks}
 \newtheorem{example}[theorem]{Example}

\newtheorem*{questions}{Questions}

\numberwithin{equation}{section}

\newcommand {\N}{\mathbb{N}} 
\newcommand {\Z}{\mathbb{Z}} 
\newcommand {\R}{\mathbb{R}} 
\newcommand {\Q}{\mathbb{Q}} 
\newcommand {\C}{\mathbb{C}} 
\newcommand{\F}{\mathbb{F}}

\newcommand{\FF}{\mathcal{F}}

\newcommand{\PP}{\mathcal{P}}

\DeclareMathOperator{\ACA}{ACA}

      \DeclareMathOperator{\SL}{SL}

\DeclareMathOperator{\Fix}{Fix}

\DeclareMathOperator{\Id}{Id}
\DeclareMathOperator{\Zer}{Z}
\DeclareMathOperator{\Ideal}{I}

\begin{document}
 \title{On algebraic cellular automata}
 \author{Tullio Ceccherini-Silberstein}
\address{Dipartimento di Ingegneria, Universit\`a del Sannio, C.so
Garibaldi 107, 82100 Benevento, Italy}
\email{tceccher@mat.uniroma1.it}
\author{Michel Coornaert}
\address{Institut de Recherche Math\'ematique Avanc\'ee,
UMR 7501,                                             Universit\'e  de Strasbourg et CNRS,
                                                 7 rue Ren\'e-Descartes,
                                               67000 Strasbourg, France}
\email{coornaert@math.unistra.fr}
\date{\today}
\begin{abstract}
We investigate some general properties of algebraic cellular automata, i.e., cellular automata over groups whose alphabets are affine algebraic sets and which are locally defined by regular maps.
When the ground field is assumed to be uncountable and algebraically closed, 
we prove that such cellular automata always have a closed image with respect to the prodiscrete topology on the space of configurations and that they are reversible as soon as they are bijective.    
\end{abstract}
\maketitle

\section{Introduction}

The goal of the present paper is to investigate some properties of algebraic cellular automata, namely cellular automata over groups whose alphabets are affine algebraic sets and whose local defining maps are regular. 
Most of the results presented here are particular cases of more general statements due to 
M. Gromov \cite{gromov-esav}.
\par 
Let us first recall some basic definitions and introduce notation. 
\par
Let $G$ be a group and let $A$ be a set.
The set $A^G = \{x \colon G \to A\}$ is called the set of \emph{configurations} over the group $G$ and the \emph{alphabet} $A$.
We equip  $A^G = \prod_{g \in G} A  $  with its \emph{prodiscrete} topology, that is, with the product topology obtained by taking the discrete topology on each factor $A$ of $A^G$.
 The continuous action of $G$ on $A^G$ defined by
$$
gx(h) = x(g^{-1}h) \qquad   \text{for all  }  g,h \in G \text{  and  } x \in A^G
$$
is called the $G$-\emph{shift} on $A^G$.
Given a configuration $x \in A^G$ and a subset $\Omega \subset G$, the element $x\vert_\Omega \in A^\Omega$ defined by $x\vert_\Omega(g) = x(g)$ for all $g \in \Omega$ is called the 
\emph{restriction} of $x$ to $\Omega$ or the \emph{pattern} of $x$ supported by $\Omega$. 
\par
A \emph{cellular automaton} over the group $G$ and the alphabet $A$ is a map 
$\tau \colon A^G \to A^G$ satisfying the following condition: 
there exist a finite subset $M \subset G$ 
and a map $\mu \colon  A^M \to A$ such that 
\begin{equation} 
\label{e;local-property}
\tau(x)(g) = \mu((g^{-1}x)\vert_M)  \quad  \text{for all } x \in A^G \text{ and } g \in G
\end{equation}
(see e.g. \cite{ca-and-groups-springer}).
 Such a set $M$ is then called a \emph{memory set} of $\tau$ and $\mu$ is called the \emph{local defining map} for $\tau$   associated with $M$. 
 \par
It easy to see from this definition that every cellular automaton $\tau \colon A^G \to A^G$ is continuous (for the prodiscrete topology on $A^G$) and $G$-equivariant 
(i.e., it satisfies $\tau(gx) = g\tau(x)$ for all $g \in G$, $x \in A^G$).
Conversely, when the alphabet $A$ is finite, it follows from the Curtis-Hedlund theorem 
(see e.g. \cite[Theorem 1.8.1]{ca-and-groups-springer}) that every continuous $G$-equivariant map from $A^G$ into itself is a cellular automaton. 
\par
Let $K$ be a field.
Recall that an \emph{affine algebraic set} over $K$  is the set of common zeroes  $\Zer(S) \subset K^m$ of a set of polynomials 
$S \subset K[t_1,\dots,t_m]$,   $m \geq 1$, and that
a map $f \colon V \to W$ between affine algebraic sets $V \subset K^m$ and $W \subset K^n$ is called \emph{regular} if $f$ is the restriction of some polynomial map $P \colon K^m \to K^n$
 (see Section \ref{sec:bacground-alg-geometry} for the required background material in affine algebraic geometry).
  One says that a cellular automaton $\tau \colon A^G \to A^G$ 
is an  \emph{algebraic cellular automaton} over the field $K$
 if the alphabet $A$ is an affine algebraic set over $K$
 and if for some (or, equivalently,  any)  memory set $M \subset G$ of $\tau$, 
 the associated local defining map $\mu \colon A^M \to A$ is regular.
\par
One says that a map $f \colon X \to Y$ between topological spaces $X$ and $Y$ has the 
\emph{closed image property} if its image $f(X)$ is closed in $Y$.
When the alphabet $A$ is a finite set, it immediately follows from the compactness of $A^G$ that every cellular automaton $\tau \colon A^G \to A^G$ has the closed image property.
 It is also true that when $A$ is a finite-dimensional vector space, then every linear cellular automaton $\tau \colon A^G \to A^G$ has the closed image property (see \cite{garden}, \cite{induction}, \cite{periodic}).
On the other hand, if $G$ is a non-periodic group (i.e., a group admitting an element of infinite order),
then, given any infinite set (resp. any infinite-dimensional vector space) $A$,  one can construct a  cellular automaton (resp. a linear cellular automaton) $\tau \colon A^G \to A^G$ whose image is not closed in $A^G$ (see \cite{periodic}).  
 \par
Our first  result on algebraic cellular automata is a particular case of  \cite[Section 4.D]{gromov-esav}:

\begin{theorem}
\label{t:closed-image}
Let $G$ be a group, $K$ an uncountable algebraically closed field, and $A$ an affine algebraic set over $K$.
Then every algebraic cellular automaton $\tau \colon A^G \to A^G$ has the closed image property with respect to the prodiscrete topology on $A^G$.
\end{theorem}

Recall that a group is called \emph{residually finite} if the intersection of its finite index subgroups is reduced to the identity element.
For example, by a theorem of Mal'cev, every finitely generated linear group is residually finite.
From Theorem \ref{t:closed-image} and the Ax-Grothendieck theorem \cite{ax-elementary}, \cite{grothendieck-4-3} on injective endomorphisms of algebraic varieties, we shall deduce the following (cf. \cite[Section 4.E']{gromov-esav}): 

\begin{corollary}
\label{c:surjunctivity-rf}
Let $G$ be a residually finite group (e.g., $G = \Z^d$), $K$ an uncountable algebraically closed field, and $A$ an affine algebraic set  over $K$.
Then every injective algebraic cellular automaton $\tau \colon A^G \to A^G$  is surjective and hence bijective.
\end{corollary}

Given a group $G$ and a  set $A$,  a cellular automaton $\tau \colon A^G \to A^G$ is called \emph{reversible} if $\tau$ is bijective and its inverse map $\tau^{-1} \colon A^G \to A^G$ is also a cellular automaton.
When the alphabet $A$ is finite, it easily follows from the compactness of $A^G$ and the 
Curtis-Hedlund theorem that every bijective cellular automaton $\tau \colon A^G \to A^G$ is reversible (see e.g. \cite[Theorem 1.10.2]{ca-and-groups-springer}).
 It is also known that when $A$ is a finite-dimensional vector space, then every bijective linear cellular automaton $\tau \colon A^G \to A^G$ is reversible (see \cite{israel}, \cite{induction}, \cite{periodic}).
On the other hand, if $G$ is a non-periodic group,
then, given any infinite set (resp. any infinite-dimensional vector space) $A$,  one can construct a bijective cellular automaton (resp. a bijective linear cellular automaton) $\tau \colon A^G \to A^G$ which is not reversible (see \cite{periodic}).  
\par
For algebraic cellular automata, we shall prove the following result:

\begin{theorem}
\label{t:reversible}
Let $G$ be a group, $K$ an uncountable algebraically closed field, and $A$ an affine algebraic set  over $K$.
Then every bijective algebraic cellular automaton $\tau \colon A^G \to A^G$ is reversible. 
\end{theorem}

The paper is organized as follows.
In Section \ref{sec:bacground-alg-geometry}, we collect the required material from affine algebraic geometry. 
In Section \ref{sec:ca}, we prove the above mentioned result  that, given   a cellular automaton whose alphabet is an affine algebraic set, the fact that the local defining map is a regular map  does not depend on the choice of the memory set. We also recall some properties of the operations of induction and restriction of a cellular automaton with respect to a subgroup of the underlying group. 
 Section \ref{sec:proj-seq} is devoted to projective sequences of constructible sets.
We prove that the projective limit of a projective sequence of nonempty  constructible sets over an uncountable algebraically closed field is never empty (Theorem \ref{t:proj-limit-constructible}).
This last result, a Mittag-Leffler-type statement,  is a key ingredient in the proof of Theorem \ref{t:closed-image} and Theorem \ref{t:reversible}.
In Section \ref{sec:closed-image}, we establish Theorem \ref{t:closed-image} 
and Corollary \ref{c:surjunctivity-rf}.
We also describe an algebraic cellular automaton over the field $\R$ with alphabet $A = \R$ and underlying group $G = \Z$ which does not have the closed image property 
(Example \ref{ex:not-closed-image}). 
The proof of Theorem \ref{t:reversible} is given in Section \ref{sec:reversibility}.
\par
The present paper grew out from numerous readings  of \cite{gromov-esav}.
Once again, we would like to express our deep gratitude to Misha Gromov for inspiration and motivation. We also thank the referee for her/his suggestions and remarks which helped us
to improve the exposition.
 
 \section{Basic affine algebraic geometry}
\label{sec:bacground-alg-geometry}

In this section, we briefly review  the material from algebraic geometry that will be needed in this paper.
For a more detailed exposition and proofs, see for example
\cite{borel-algebraic}, \cite{bump}, \cite{hartshorne}, \cite{kraft}, \cite{milne} or \cite{reid}.

\subsection{Affine algebraic sets}
Let $K$ be a field and let $m$ be a positive  integer.
Let $K[t_1,\dots,t_m]$ denote the $K$-algebra of  polynomials  in $m$ indeterminates with coefficients in $K$. 
Given a subset $S \subset K[t_1,\dots,t_m]$, we denote by  $\Zer(S)$ the subset of $K^m$ consisting  of the common zeroes for the polynomials in $S$, that is,
$$
\Zer(S) = \{a = (a_1,\dots,a_m) \in K^m : P(a) = 0 \quad \text{for all} \quad P \in S\}.
$$
When $S = \{P_1,\dots,P_s\}$ is a finite set, we shall also write
$\Zer(P_1,\dots,P_s)$ instead of $\Zer(S)$.
\par
One says that a subset $A \subset K^m$ is an \emph{algebraic subset} of $K^m$ if there exists a subset $S \subset K[t_1,\dots,t_m]$ such that $A = \Zer(S)$.
An algebraic subset $A \subset K^m$ is also called an \emph{affine algebraic set} or an 
\emph{affine algebraic variety} (for some authors an affine algebraic variety must be irreducible in the sense of Subsection \ref{subsec:irr} below).

 The intersection of any family of algebraic subsets of $K^m$
as well as the union of any finite family of algebraic subsets of $K^m$
are also algebraic subsets of $K^m$.
 It follows that the algebraic subsets of $K^m$ are the closed subsets of a topology on $K^m$.
This topology is called the \emph{Zariski topology} on $K^m$.
In the sequel, unless another topology is specified, it will be understood that the topology on 
$K^m$ (resp. on any subset of $K^m$) is the Zariski topology (resp. the topology induced by the Zariski topology of $K^m$).
\par
 Given an arbitrary subset $\Sigma$ of $K^m$, the set
\begin{equation}
\label{def:ideal-subset}
\Ideal(\Sigma) = \{P \in K[t_1,\dots,t_m] : P(a) = 0 \quad \text{for all} \quad a \in \Sigma \} \subset K[t_1,\dots,t_m],
\end{equation}
consisting  of all polynomials which are identically $0$  on $\Sigma$,
 is an ideal of $K[t_1,\dots,t_m]$.
 One has
 \begin{equation}
\label{e:Sigma-included-Z}
\Sigma \subset \Zer(\Ideal(\Sigma))
\end{equation}
for every subset $\Sigma \subset K^m$ and equality holds if and only if $\Sigma$ is an algebraic subset of $K^m$.
 \par
Let $A \subset K^m$ be an algebraic subset.
The quotient ring $K[t_1,\dots,t_m]/\Ideal(A)$ is denoted by $K[A]$ and is called the 
\emph{coordinate ring} of $A$. In fact, $K[A]$ inherits from
$K[t_1,\dots,t_m]$ a structure of a  $K$-algebra.
It can be viewed as a sub-$K$-algebra of the $K$-algebra formed by all $K$-valued maps on $A$.
\par
  As the ring $K[t_1,\dots,t_m]$ is Noetherian,
one can find finitely many polynomials
$P_1,\dots,P_r \in K[t_1,\dots,t_m]$ generating the ideal $\Ideal(A)$.
We then have 
$$
A = \Zer(P_1,\dots,P_r) = \bigcap_{1 \leq i \leq r} \Zer(P_i).
$$
As the closed subsets of $A$ are the algebraic subsets of $K^m$ which are contained in $A$, a subset $B \subset A$ is closed if and only if there exist finitely many polynomials $Q_1,\dots,Q_s \in K[t_1,\dots,t_m]$ such that
$$
B = A \cap \Zer(Q_1, \dots,Q_s)= \Zer(P_1,\dots,P_r,Q_1,\dots,Q_s).
$$
By taking complements, we deduce that $U \subset A$ is open if and only if there exist finitely many polynomials $Q_1,\dots,Q_s \in K[t_1,\dots,t_m]$ such that
\begin{equation}
\label{e:open-union-special}
 U = U_{Q_1} \cup \dots \cup U_{Q_s},
\end{equation}
where we set, for $Q \in K[t_1,\dots,t_m]$,
$$
U_Q = A \setminus \Zer(Q) = \{ a \in A : Q(a) \not= 0\}.
$$
In the case when $U = U_Q$ for some $Q \in K[t_1,\dots,t_m]$,
one says that $U$ is a \emph{special open subset} of $A$.
Thus, it follows from \eqref{e:open-union-special}
 that every open subset of $A$ is a finite union of special open subsets.
\par
The fact that $K[t_1,\dots,t_m]$ is Noetherian also implies that the algebraic subsets of 
$K^m$ satisfy the \emph{descending chain condition}. This means that every sequence 
 $(A_n)_{n \geq 1}$ of algebraic subsets of $K^m$  with 
 $A_{n + 1} \subset A_n$ for all $n \geq 1$ eventually stabilizes. 
\par
In the case when $K$ is algebraically closed, the map $A \mapsto \Ideal(A)$ is an inclusion-reversing bijection from the set of algebraic subsets of $K^m$ onto the set of radical ideals in 
$K[t_1,\dots,t_m]$.

 \subsection{Irreducible components}
\label{subsec:irr}
An algebraic subset $A \subset K^m$ is called \emph{irreducible}
if it is nonempty and it cannot be expressed in the form $A = A_1 \cup A_2$ where $A_1$ and $A_2$ are proper closed subsets of $A$.
This is equivalent to the fact that the ideal $\Ideal(A)$ is prime in the ring $K[t_1,\dots,t_m]$.
\par  
Every algebraic subset $A \subset K^m$ can be expressed as a finite union
$A = A_1 \cup A_2 \cup \dots \cup A_s$,
where $A_1$, $A_2$, \dots , $A_s$ are irreducible algebraic subsets of $K^m$ such that there are no $i \not= j$ with  $A_i  \subset A_j$.
  Such a decomposition is unique up to reordering. The algebraic subsets $A_1$, $A_2$, \dots, $A_s$ are the maximal closed irreducible subsets of $A$ and are called the
\emph{irreducible components} of $A$.
 \par
 When $A \subset K^m$ is an irreducible algebraic subset, the field of fractions of its coordinate ring $K[A]$ is called the \emph{function field} of $A$ and is denoted by $K(A)$.

\subsection{Regular maps}
Let $m$ and $n$ be positive integers.
One says that a map $F \colon K^m \to K^n$ is \emph{polynomial} if there exist polynomials $P_1,\dots,P_n \in K[t_1,\dots,t_m]$ such that $F(a) = (P_1(a),\dots,P_n(a))$ for all $a = (a_1,\dots,a_m)  \in K^m$.
Let $A \subset K^m$ and $B \subset K^n$ be algebraic subsets. One says that a map 
$f \colon A \to B$ is \emph{regular} if $f$ is the restriction of some polynomial map  
$F \colon K^m \to K^n$.
Every regular map $f \colon A \to B$ is continuous for the Zariski topology and induces a $K$-algebra homomorphism
$f^* \colon K[B] \to K[A]$ given by $f^*(\varphi) = \varphi \circ f$ for all $\varphi \in K[B]$.
\par 
The identity map $\Id_A \colon A \to A$ is regular
for any algebraic subset $A \subset K^m$.
Moreover,  if the maps $f \colon A \to B$ and $g \colon B \to C$ are regular, where $A \subset K^m$, $B \subset K^n$, and $C \subset K^p$ are algebraic subsets, then the composite map 
$g \circ f \colon A \to C$ is also regular. 
It follows that the algebraic subsets of $K^m$, $m = 1,2,\dots$, are the objects of a category whose morphisms are the regular maps between them. This category is called the 
 \emph{category of affine algebraic sets} over $K$. 
\par
The category of affine algebraic sets over $K$ admits finite direct products.
Indeed, if $A \subset K^m$ and $B \subset K^n$ are algebraic subsets,
then  the Cartesian product $A \times B$ is an algebraic subset of  $ K^m \times K^n = K^{m + n}$,
and  the two projection maps $\pi_A \colon A \times B \to A$ and $\pi_B \colon A \times B \to B$ are regular. 
One easily checks that the triple $(A \times B, \pi_A,\pi_B)$ is the direct product of $A$ and $B$ in the category of affine algebraic sets over $K$.

\subsection{Constructible sets and Chevalley's theorem}
\label{subsec:construcible}
Let $A$ be an affine algebraic set over a field $K$.
One says that a subset $L \subset A$ is \emph{locally closed} in $A$ if there exist an open subset 
$U \subset A$   and a closed subset $V \subset A$ such that $L = U \cap V$.
This is equivalent to  $L$ being open in its closure $\overline{L} \subset A$.
 \par
   One says that a subset $C \subset A$ is \emph{constructible} if $C$ is a finite union of locally closed subsets of $X$.
  The set of constructible subsets of $A$ is closed under finite unions, finite intersections,
and taking complements in $A$.
 It follows that the set of constructible subsets of $A$ is a Boolean subalgebra of the Boolean algebra $\PP(A)$ formed by  all subsets of $A$.
It is the smallest Boolean subalgebra of $\PP(A)$ containing all closed subsets of $A$.
We shall use the following elementary result 
(see for example \cite[AG Section 1.3]{borel-algebraic}):

 \begin{proposition}
 \label{p:constructible-contains-open}
Let $A$ be an affine algebraic set over a field $K$ and suppose  that $C$ is a constructible subset of $A$.
 Then there is an open dense subset $U$ of $\overline{C}$ such that $U \subset C$.
\end{proposition}

We shall also use the following theorem due to C. Chevalley (see for example 
\cite[AG Section 10.2]{borel-algebraic}):

\begin{theorem}[Chevalley]
\label{t:chevalley}
Let $K$ be an algebraically closed field. Let $A $ and $B$ be affine algebraic sets over $K$,  and let $f \colon A \to B$ be a regular map.
 Then every constructible subset $C \subset A$  has a constructible image $f(C) \subset B$. 
In particular, $f(A)$ is a constructible subset of $B$.
\end{theorem}

\subsection{Dimension}
In this subsection, the field $K$ is  assumed to be algebraically closed.
Let $A$ be an affine algebraic set over $K$.
  The \emph{ dimension} $\dim(A)$ of $A$  is defined as being the greatest  integer  $n \geq 0$ such that there exists a strictly increasing chain $(F_i)_{0 \leq i \leq n}$
of length $n$ consisting of closed irreducible subsets of $A$
(by convention, the  dimension of the empty set   is $-\infty$).
\par
 One has $\dim(K^m) = m$.
 If $A = A_1 \cup \dots \cup A_s$ is the decomposition of $A$ into irreducible components, then 
$\dim(A) = \max_{1 \leq i \leq s} \dim(A_i)$. 
If $B$ is a closed subset of $A$, then one always has $\dim(B) \leq \dim(A)$.
Moreover, if $A$ is irreducible and $B$ is a closed subset of $A$ with $B \not= A$ then one has 
$\dim(B) < \dim(A)$.
If $f \colon A \to A'$ is a surjective regular map between affine algebraic sets, then 
$\dim(A') \leq \dim(A)$.
\par
The dimension of an affine algebraic set $A$ is equal to the Krull dimension of its coordinate ring $K[A]$, i.e., to the maximal length of a strictly increasing chain of prime ideals of $K[A]$.
If in addition  $A$ is  irreducible,   then $\dim(A)$ is also equal to the transcendence degree of its function field $K(A)$ over $K$.
\par
Let $A$ and $B$ be irreducible affine algebraic sets over $K$.
Let $f \colon A \to B$ be a regular map and let 
$f^* \colon K[B] \to K[A]$ denote the induced ring homomorphism.
 One says that $f$ is a \emph{finite morphism} if $K[A]$ is finitely generated as a $f^*(K[B])$-module.
Every finite morphism $f \colon A \to B$ between irreducible affine algebraic sets is closed, i.e., such that the image of any closed subset of $A$ is closed in $B$ (see for example \cite[Proposition 8.7]{milne}).
\par 
We shall use the following result (see for example \cite[Theorem 8.12]{milne}) which can be deduced from Emmy Noether's normalization lemma:

\begin{theorem}
\label{t:finite-morphism}
Let $K$ be an algebraically closed field and let  $A$ be an irreducible affine algebraic set over $K$ such that    $ \dim(A) = d$.
Then there exists a surjective finite morphism $f \colon A \to K^d$.
\end{theorem}  

  \subsection{The Ax-Grothendieck theorem}
\label{subsec:ax}

In the proof of Corollary \ref{c:surjunctivity-rf},
we shall use the following result:

\begin{theorem}[Ax-Grothendieck]
\label{t:ax-theorem}
Let $K$ be an algebraically closed field  and let $A$ be an affine algebraic set over $K$.
Then every injective regular map $f \colon A \to A$ is surjective and hence bijective.
\end{theorem}

Theorem \ref{t:ax-theorem} was established independently by J. Ax \cite[Theorem C]{ax-elementary}, \cite{ax-injective}  
and by A. Grothendieck 
\cite[Proposition 10.4.11]{grothendieck-4-3} in the more general setting of injective endomorphisms of schemes of finite type.  
The proof of Ax is model-theoretic. 
A cohomological proof of the Ax-Grothendieck theorem for algebraic varieties was  given by 
A. Borel in \cite{borel-basel}.
An elementary proof of Theorem \ref{t:ax-theorem} may be found in \cite{kang} 
(see also \cite{van-den-essen}). 
The Ax-Grothendieck theorem is also discussed in \cite{gromov-esav} and \cite{serre-how-to-use}.

\begin{remarks}
\label{rem:ax}
(a)
Theorem \ref{t:ax-theorem} becomes false if the hypothesis that $K$ is algebraically closed is removed.
\par
In characteristic $0$, it suffices to consider  
  the injective polynomial map $f \colon \Q \to \Q$ 
  defined by $f(x) = x^3$, which  is not surjective since $2 \notin f(\Q)$.
 \par
 In positive characteristic, examples of injective but not surjective polynomial self-mappings of fields may be obtained as follows. 
  Let $K$ be a field of characteristic $p > 0$ and consider  the polynomial map $f \colon K \to K$ defined by $f(x) = x^p$ for all $x \in K$.
One clearly has $f(1_K) = 1_K$ and $f(xy) = f(x)f(y)$ for all $x,y \in K$.
Moreover, the binomial formula applied to $(x + y)^p$ shows that one also has
$f(x + y) = f(x) + f(y)$ for all $x,y \in K$, since $p$ divides ${p \choose k}$ 
for all $1 \leq k \leq p - 1$.
It follows that $f$ is an endomorphism of the field $K$. In particular, $f$ is injective.
This endomorphism  is called the \emph{Frobenius endomorphism} of $K$. 
The Frobenius endomorphism may fail to be surjective.
For instance, if $k$ is any field of characteristic $p > 0$ (e.g., $k = \Z/p\Z$)
and $K = k(t)$ denotes  the field of rational functions with coefficients in  $k$ in one indeterminate $t$, then the Frobenius endomorphism $f \colon K \to K$   is not surjective since there is no 
$R \in K$ such that $t = R^p$. 
\par
(b) When $A$ is an affine algebraic set over an algebraically closed field $K$ of characteristic $0$, it is known that the inverse map of any bijective regular map $f \colon A \to A$ is also regular
(see \cite[Proposition 17.9.6]{grothendieck-ega-1667}, \cite{cynk}).
\par
(c) When $K$ is an algebraically closed field of characteristic $p > 0$ and $A$ is an affine algebraic set over $K$, the inverse map of a bijective regular map $f \colon A \to A$
need not to be regular. For example, the inverse map of the Frobenius automorphism 
$f \colon K \to K$ is not regular since there is no polynomial $P in K[t]$ such that
$P(x)^p = x$ for all $x \in K$.
\par
(d) 
It is known \cite{kurdyka} that every injective regular map $f \colon A \to A$, where $A$ is a real affine algebraic set, is bijective. 
However, its inverse need not to be regular.
For example, the inverse of the bijective polynomial map $f \colon \R \to \R$ defined by $f(x) = x^3$ is the map
$x \mapsto \sqrt[3]{x}$ which is not polynomial.
\par
(e)
The inverse map    of a bijective regular map between distinct algebraic subsets may fail to be regular even if the ground field is algebraically closed and of characteristic $0$.
For example, the map $f \colon \C \to \Zer(t_1^2 - t_2^3) \subset \C^2$ given by $f(t) = (t^3,t^2)$ is bijective and regular but its inverse map
 is not regular. 
  Otherwise, this would imply the existence of a polynomial
$P \in \C[t_1,t_2]$ such that $P(z^3,z^2) = z$ for all $z \in \C$.
This is impossible since, for any $P \in \C[t_1,t_2]$, the expression  $P(z^3,z^2)$ is polynomial in $z$ with each non-constant monomial of degree at least $2$. 
\end{remarks}

\section{Cellular automata}
\label{sec:ca}

\subsection{Algebraic cellular automata}

Let $G$ be a group and let $A$ be a set.
Given a cellular automaton $\tau \colon A^G \to A^G$ and a memory set $M \subset G$ for $\tau$, we denote by
$\mu_M \colon A^M \to A$ the local defining map for $\tau$ associated with $M$.
Observe that $\mu_M$ is entirely determined by $\tau$ and $M$ since, for all $y \in A^M$, we have
\begin{equation}
\label{e:mu-M}
\mu_M(y) = \tau(x)(1_G),
\end{equation}
where $x \in A^G$ is any configuration satisfying $x\vert_M = y$.
\par
We recall (see for example \cite[Section 1.5]{ca-and-groups-springer}) that every cellular
automaton $\tau \colon A^G \to A^G$ admits a unique memory set
$M_0 \subset G$ of minimal cardinality and that in addition a
subset $M \subset G$ is a memory set for $\tau$ if and only
if $M_0 \subset M$. Such a memory set is called the \emph{minimal}
memory set of $\tau$.

\begin{proposition}
\label{p:carac-aca}
Let $G$ be a group and let $A$ be an affine algebraic set over a field $K$.
Let $\tau \colon A^G \to A^G$ be a cellular automaton.
Then the following conditions are equivalent:
\begin{enumerate}[\rm (a)]
\item
there exists a memory set $M$ of $\tau$ such that 
the associated local defining map  $\mu_M \colon A^M \to A$ is regular;
\item
for any memory set $M$ of $\tau$,
the associated local defining map  $\mu_M \colon A^M \to A$ is regular.
\end{enumerate}
\end{proposition}

\begin{proof}
Suppose that the local defining map $\mu_M \colon A^M \to A$ is regular for some memory set 
$M$ of $\tau$.
Let $M'$ be another memory set of $\tau$ and let us show that the associated local defining map 
$\mu_{M'} \colon A^{M'} \to A$ is also regular.
Consider the minimal memory set $M_0$ of $\tau$ 
 and fix an arbitrary point $a_0 \in A$.
 We have $M_0 \subset M$ and
$\mu_{M_0} = \mu_M \circ \iota$, where $\iota \colon A^{M_0} \to A^M$ is the
embedding defined by 
$$
\iota(y)(g) =
\begin{cases}
y(g) & \text{ if  } g \in M_0, \\
a_0 & \text{ if  } g \in M \setminus M_0,
\end{cases}
$$
for all $y \in A^{M_0}$. It follows that the map $\mu_{M_0}$ is regular.
On the other hand, we have $M_0 \subset M'$ and
$\mu_{M'} = \mu_{M_0} \circ \pi$, where $\pi \colon A^{M'} \to A^{M_0}$ is the projection map.
We deduce that  $\mu_{M'}$ is a regular map.
\end{proof}

Given a field $K$, we say that a cellular automaton $\tau \colon A^G \to A^G$ is  an 
\emph{algebraic cellular automaton} over $K$ if $A$ is an algebraic set over $K$ and $\tau$ satisfies one of the equivalent conditions of Proposition \ref{p:carac-aca}.

\begin{example}{\rm 
Every cellular automaton with finite alphabet $A$ may be regarded as an algebraic cellular automaton. Indeed, it suffices to embed $A$ as a subset of some field $K$ and then observe that, if $M$ is a finite set, any map $\mu \colon A^M \to A$ is the restriction of some polynomial map $P \colon K^M \to K$ (which can be made
explicit by using Lagrange interpolation formula).} 
\end{example}

\begin{example}
\label{ex:premier-exemple}
{\rm Let $K$ be a field, $A$ an affine algebraic set over $K$, and $f \colon A \to A$ a regular map. Let $G$ be a group and fix an element $g_0 \in G$. Then the map $\tau \colon A^G \to A^G$, defined by
\begin{equation}
\label{e:premier-exemple}
\tau(x)(g) = f(x(gg_0))
\end{equation}
for all $x \in A^G$ and $g \in G$, is an algebraic cellular automaton with memory
set $\{g_0\}$ and local defining map $f$ (we have identified $A^{\{g_0\}}$ with $A$).
Note that $\tau$ is injective (resp. surjective) if and only if $f$ is injective (resp.
surjective).}
\end{example}

\begin{example}
{\rm Let $K$ be a field. Let $A$ be an affine algebraic group over $K$, i.e., an algebraic set over $K$ equipped with a group structure such that both the group multiplication and the inverse operation are given by regular maps (for example $A = \SL_N(K)$).
Then the map $\tau \colon A^\Z \to A^\Z$, defined by 
$$
\tau(x)(n) = (x(n + 1))^{-1} x(n)
$$
for all $x \in A^\Z$ and $n \in \Z$, is an algebraic cellular automaton.
Note that $\tau$ is surjective and that $\tau$ is not injective unless $A$ is reduced to the identity element.}
\end{example}

\begin{example}
{\rm Let $K$ be an algebraically closed field.
Then a map $\tau \colon K^\Z \to K^\Z$ is an injective algebraic cellular automaton
if and only if there exists $m_0 \in \Z$ and $\alpha,\beta \in K$, $\alpha \neq 0$, such that
\begin{equation}
\label{e:memoire}
\tau(x)(n) = \alpha x(n+m_0) + \beta
\end{equation}
for all $x \in K^\Z$ and $n \in \Z$. Indeed, if $\tau$ is of the form
\eqref{e:memoire} then we are exactly in the situation described in
Example \ref{ex:premier-exemple} with $G = \Z$, $g_0 = m_0$, $A = K$ and $f(z) = \alpha z + \beta$.
As $f \colon K \to K$ is an injective affine map, it follows that $\tau$ is an injective algebraic cellular automaton.
\par
Conversely, suppose that $\tau \colon K^\Z \to K^\Z$ is an injective
algebraic cellular automaton. Let $M \subset \Z$ be the minimal memory set for $\tau$ and denote by $\mu \colon K^M \to K$ the corresponding local defining map. Suppose by contradiction that $M$ has cardinality $N \geq 2$ and let $m_1 < m_2 < \cdots < m_N$ denote its elements. We construct a configuration $y \in K^\Z$ as follows. First, for $k \in K$, let $x_k \in K^\Z$ denote the constant configuration defined by $x_k(n) = k$ for all $n \in \Z$. Thus if $c = \tau(x_0)(0) \in K$ we have $\tau(x_0) = x_c$. Then, we choose arbitrary values $y(m_1), y(m_1+1), \ldots, y(m_N-1) \in K$ such that $y(m_1) \neq 0$.  Since $K$ is algebraically closed, we can find $b \in K$ such that $\mu(y(m_1), y(m_2), \cdots, y(m_{N-1}),b) = c$. We then set $y(m_N) = b$. Similarly, we can find $b' \in K$ such that $\mu(y(m_1+1), y(m_2+1), \cdots, y(m_{N-1}+1),b') = c$. We then set $y(m_N+1) = b'$. Continuing this way, we define all the values $y(n)$ for $n \geq m_1$. Symmetrically, we can find $b'' \in K$ such that $\mu(b'',y(m_2-1), y(m_3-1), \cdots, y(m_N-1)) = c$. We then set $y(m_1-1) = b''$. Continuing this way,  all the values $y(n)$ with $n \leq m_1-1$ are also defined. By construction, we have $y \neq x_0$ (since $y(m_1) \neq 0 = x_0(m_1)$). Moreover $\tau(y) = x_c = \tau(x_0)$. This contradicts the injectivity of $\tau$. We have shown that $\vert M \vert = 1$. In this case, the injectivity of $\tau$ is equivalent to the injectivity of $\mu \colon K \to K$ so that the local defining map is necessarily of the form
$\mu(z) = \alpha z + \beta$ for suitable $\alpha, \beta \in K$ with $\alpha \neq 0$. This
shows \eqref{e:memoire}, where $M =\{m_0\}$. Note that $\tau$ is in fact bijective with inverse map $\tau^{-1} \colon K^\Z \to K^\Z$ given by  $\tau^{-1}(x)(n) = \alpha^{-1}x(n-m_0) - \alpha^{-1}\beta$ for all $x \in K^\Z$ and $n \in \Z$. Thus $\tau^{-1}$
is an algebraic cellular autmaton as well (with memory set $\{-m_0\}$).
\par
Consequently, when $K$ is algebraically closed, every injective cellular automaton
$\tau \colon K^\Z \to K^\Z$ is bijective and it inverse map $\tau^{-1} \colon K^\Z \to K^\Z$
is also an algebraic cellular automaton.}
\end{example}

\begin{example}
{\rm Let $K$ be a field. Let  $M = \{1,2,\ldots, m\} \subset \Z$ and $A = K^m$, where $m \geq 2$.
Given a configuration $x \in A^\Z$ we write $x(n) = (x_1(n), x_2(n), \ldots, x_m(n)) \in K^m$ for all $n \in \Z$. For $i=1,2,\ldots, m$, we arbitrarily choose $\alpha_i \in K \setminus \{0\}$ and  $P_i \in K[t_1,t_2, \ldots, t_{i-1}]$.
Consider the map $\tau \colon A^\Z \to A^\Z$ defined by
$\tau(x) = y$ where
\begin{equation}
\label{e:def-tau}
\left\{
\begin{split}
y_1(n) & = \alpha_1 x_1(n+1) + P_1\\
y_2(n) & = \alpha_2 x_2(n+2) + P_2(x_1(n+1))\\
y_3(n) & = \alpha_3 x_3(n+3) + P_3(x_1(n+1),x_2(n+2))\\
 \cdots & \\
y_m(n) & = \alpha_m x_m(n+m)+ P_m(x_1(n+1),x_2(n+2), \ldots, x_{m-1}(n+m-1))
\end{split}
\right.
\end{equation}
for all $x \in A^\Z$ and $n \in \Z$.
Then $\tau$ is an algebraic cellular automaton with memory set $M$.
\par
From \eqref{e:def-tau} we immediately deduce that $\tau$ is bijective with inverse map
$\tau^{-1} \colon A^\Z \to A^\Z$ given by
$\tau^{-1}(y) = x$ where
\begin{equation}
\label{e:def-tau-inv}
\left\{
\begin{split}
x_1(n) & = \alpha_1^{-1}y_1(n-1) + Q_1\\
x_2(n) & = \alpha_2^{-1}y_2(n-2) + Q_2(y_1(n-2))\\
x_3(n) & = \alpha_3^{-1}y_3(n-3) + Q_3(y_1(n-3),y_2(n-3))\\
 \cdots & \cdots\\
x_m(n) & = \alpha_m^{-1}y_m(n-m) + Q_m(y_1(n-m),y_2(n-m), \ldots, y_{m-1}(n-m))
\end{split}
\right.
\end{equation}
for all $y \in A^\Z$ and $n \in \Z$, where the polynomials $Q_i \in K[t_1,t_2, \ldots, t_{i-1}]$ are recursively given by
\[
\left\{
\begin{split}
& Q_1 = - \alpha_1^{-1} P_1\\
& Q_2(t_1)  = - \alpha_2^{-1} P_2(\alpha_1^{-1}t_1 + Q_1)\\
& Q_3(t_1,t_2) = - \alpha_3^{-1} P_3(\alpha_1^{-1}t_1 + Q_1,\alpha_2^{-1}t_2 + Q_2(t_1))\\
& \cdots\\
& Q_m(t_1,t_2, \ldots, t_{m-1}) = - \alpha_m^{-1} P_m(\alpha_1^{-1}t_1 + Q_1,\alpha_2^{-1}t_2 + Q_2(t_1), \ldots\\
& \hspace{4.5cm} \ldots, \alpha_{m-1}^{-1}t_{m-1} + Q_{m-1}(t_1,t_2, \ldots, t_{m-2})).
\end{split}
\right.
\]
This shows that $\tau^{-1}$ is an algebraic cellular automaton as well (with memory set
$\{-m,-m+1, \ldots, -2,-1\}$).}
\end{example}

\begin{remark}
Let $G$ be a group, $K$ a field, $A$ an affine algebraic set over $K$, and $\tau \colon A^G \to A^G$ and $\sigma \colon A^G \to A^G$ two algebraic cellular automata. Then the composition $\sigma \circ \tau \colon A^G \to A^G$ is again an algebraic cellular automaton.
The fact that $\sigma \circ \tau$ is a cellular automaton is well known (see, for instance, \cite[Proposition 1.4.9]{ca-and-groups-springer}). To see that $\sigma \circ \tau$ is algebraic, we recall from \cite[Remark 1.4.10]{ca-and-groups-springer} the following facts. If $T$ (resp. $S$) is a memory set for $\tau$ (resp. $\sigma$) and $\mu\colon A^T \to A$ and $\nu \colon A^S \to A$ are the corresponding local defining maps, then $ST = \{st: s \in S, t \in T\}$ is a memory set for $\sigma \circ \tau$ and the corresponding local defining map can be described as follows.
For $y \in A^{ST}$ and $s \in S$, define $y_s \in A^T$ by setting $y_s(t) = y(st)$
for all $t \in T$. Also denote by $\overline{y} \in A^S$ the map defined by $\overline{y}(s) = \nu(y_s)$ for all $s \in S$. Then the local defining map for
$\sigma \circ \tau$ is the map $\kappa \colon A^{ST} \to A$ given by
$\kappa(y) = \mu(\overline{y})$ for all $y \in A^{ST}$. Now, since the maps $\nu$ and $y \mapsto y_s$, $s \in S$, are regular, we have that the map $y \to \overline{y}$ is also regular. Composing the latter with the regular map $\mu$ we obtain $\kappa$ which is therefore regular as well.
\par
Since the identity map $\Id_{A^G}\colon A^G \to A^G$ is an algebraic cellular automaton,
we have that the set $\ACA(G;A)$ consisting of all algebraic cellular automata $\tau \colon A^G \to A^G$ is a monoid for the composition of maps.
\end{remark}

\subsection{Induction and restriction}
Let $G$ be a group, $A$ a set, and $H$ a subgroup of $G$.
\par
Suppose that a cellular automaton $\tau \colon A^G \to A^G$ admits a memory set $M$ such 
that $M \subset H$. Let $\mu \colon A^M \to A$ denote the associated local defining map.
Then the map $\tau_H \colon A^H \to A^H$ defined by
$$
\tau_H(y)(h) = \mu((h^{-1}y)\vert_M)
\quad \text{  for all  } y \in A^H \text{ and } h \in H,
$$
is a cellular automaton over the group $H$ and the alphabet $A$, with memory set $M$ and local defining map $\mu$.
One says that $\tau_H$ is the cellular automaton obtained by \emph{restriction} of  
$\tau$ to $H$.
\par
Conversely, suppose that  $\sigma \colon A^H \to A^H$  is a cellular automaton with memory set $N \subset H$ and local defining map
$\nu \colon A^N \to A$. Then the map $\sigma^G \colon A^G \to A^G$ defined by
$$
\sigma^G(x)(g) = \nu((g^{-1}x)\vert_N)
\quad \text{ for all } x \in A^G \text{  and  } g \in G,
$$
is a cellular automaton over the group $G$ and the alphabet $A$, with memory set $N$ and local defining map $\nu$.
One says that $\sigma^G$ is the cellular automaton obtained by \emph{induction} of $\sigma$ to $G$.
\par
It immediately follows from their definitions that induction and restriction are operations one inverse to the other in the sense that one has $(\tau_H)^G = \tau$ 
and $(\sigma^G)_H = \sigma$ for every cellular automaton  $\tau \colon A^G \to A^G$ over $G$ admitting a memory set contained in $H$ and every cellular automaton $\sigma \colon A^H \to A^H$ over $H$.
We shall use the following result:

\begin{theorem}
\label{t:induction}
Let $G$ be a group, $A$ a set, and $H$ a subgroup of $G$. 
Suppose that  $\tau \colon A^G \to A^G$ is a cellular automaton over $G$ admitting
a memory set contained in $H$
and let $\tau_H \colon A^H \to A^H$ denote the cellular automaton over $H$ obtained by restriction. 
Then the following holds:
\begin{enumerate}[{\rm (i)}]
 \item 
$\tau$ is bijective if and only if $\tau_H$ is bijective;
\item 
$\tau$ is reversible if and only if $\tau_H$ is reversible;
 \item 
$\tau(A^G)$ is closed in $A^G$ for the prodiscrete topology  if and only if $\tau_H(A^H)$ is closed in $A^H$ 
for the prodiscrete topology; 
\item
when $A$ is an affine algebraic set over a field $K$, then $\tau$ is algebraic if and only if $\tau_H$ is algebraic.
\end{enumerate}
\end{theorem}

\begin{proof}
Assertions (i), (ii), and (iii) are established in \cite[Theorem 1.2]{induction}.
Assertion (iv) immediately follows from the definition of an algebraic cellular automaton since $\tau$ and $\tau_H$ admit a common local defining map. 
\end{proof}

\section{Projective sequences of constructible sets}
\label{sec:proj-seq}

Let $\N$ denote the set of nonnegative integers.
A  \emph{projective sequence} of sets is a sequence $(X_n)_{n \in \N}$ of sets equipped with maps $f_{nm} \colon X_m \to X_n$, defined for all $n,m \in \N$ with $ m \geq n$,
satisfying the following conditions:

\begin{enumerate}[(PS-1)]
\item
$f_{n n}$ is the identity map on $X_n$ for all $n \in \N$;
\item
$f_{n k} = f_{n m} \circ f_{m k}$ for all $n,m,k \in \N$ such that $k \geq m \geq n$.
\end{enumerate}
We shall denote such a projective sequence by $(X_n,f_{n m})$ or simply by $(X_n)$.
Observe that the projective sequence $(X_n,f_{n m}$) is entirely determined by the maps
$$
g_n = f_{n,n+1} \quad (n \in \N)
$$
since
\begin{equation}
\label{e:gives-fnm-gn}
f_{n m} = g_n \circ g_{n + 1} \circ \ldots \circ g_{m - 1} 
\end{equation}
for all $m > n$. Conversely, if we are given a sequence of maps $g_n \colon X_{n + 1} \to X_n$, $n \in \N$, then there is a unique projective sequence $(X_n, f_{n m})$ satisfying \ref{e:gives-fnm-gn}.
\par
Let $(X_n,f_{n m})$ be a projective sequence of sets.
 The \emph{projective limit} $X = \varprojlim X_n$ of the projective sequence $(X_n,f_{n m})$ is the subset  $X \subset \prod_{n \in \N} X_n$ consisting of the sequences $x = (x_n)_{n \in \N}$ satisfying 
$x_n = f_{n m}(x_m)$ for all $n,m \in \N$ such that $m \geq n$.
Note that there is a canonical map $\pi_n \colon X \to X_n$ sending $x$ to $x_n$ and that one has
$\pi_n = f_{n m} \circ \pi_m$ for all $m,n \in \N$ with $m \geq n$. 
\par
Property (PS-2) implies that, for each $n \in \N$, the sequence of sets $f_{nm}(X_m)$, $m \geq n$, is non-increasing. Let us set, for each $n \in \N$, 
$$
X_n' = \bigcap_{m \geq n} f_{nm}(X_m).
$$
The set $X_n'$ is called the set of \emph{universal elements} in $X_n$ 
(cf. \cite{grothendieck-ega-3}).
Observe that $f_{n m}(X_m') \subset X_n'$ for all $m \geq n$.
Thus, the map $f_{n m}$ induces by restriction a map $f_{n m}' \colon X_m' \to X_n'$ 
for all $m  \geq n$. 
Then $(X_n',f_{n m}')$ is a projective sequence which is called the
\emph{universal projective sequence} associated with the projective sequence $(X_n,f_{n m})$.
It is clear that the projective sequences $(X_n,f_{n m})$ and $(X_n',f_{n m}')$ have the same projective limit.
\par
The following result belongs to the prosperous family of Mittag-Leffler-type statements
(see e.g. \cite[TG II. Section 5]{bourbaki-top-gen}, \cite[Section I.3]{grothendieck-ega-3}, \cite[Section 3]{periodic}).

\begin{proposition}
\label{p:universal}
Let $(X_n,f_{n m})$ be a projective sequence of sets and let $(X_n',f_{n m}')$ denote the associated universal projective sequence of sets.
Let $X = \varprojlim X_n = \varprojlim X_n'$ denote their common projective limit.
Suppose that all maps $f_{n m}' \colon X_m' \to X_n'$, $m, n \in \N$ and $m \geq n$ are surjective.
Then all canonical maps $\pi_m' \colon X \to X_n'$, $m \in \N$, are surjective.
In particular, if   $X_m' \not= \varnothing$ for all $m \in \N$, then one has $X \not= \varnothing$. 
\end{proposition}

\begin{proof}
Let $x_m' \in X_m'$.
As the maps $f_{k,k + 1}'$, $k \geq m$, are surjective, we can construct by induction a sequence 
$(x_k')_{k \geq m}$ such that $x_k' = f_{k,k + 1}'(x_{k + 1}')$ for all $k \geq m$.
Then the sequence $(x_n)_{n \in \N}$, where $x_n = x_n'$ if $n \geq m$ and
$x_n = f_{n m}(x_m')$ if $n < m$, is in $X$ and satisfies $x_m' = \pi_m'(x)$. This shows that 
$\pi_m'$ is surjective.
 \end{proof}

\begin{remark}
Observe that, for the maps $f_{n m}'$, $m \geq n$, to be surjective,  it suffices that  the maps 
$f_{n,n+1}'$ are surjective. Also, for the sets $X_n'$ to be nonempty, $n \in \N$,  it suffices that the set $X_0'$ is nonempty. 
\end{remark}

Let $K$ be a field.
We say that a projective sequence $(X_n,f_{n m})$ is a \emph{projective sequence of constructible sets} over $K$ if there is a projective sequence $(A_n,F_{n m})$ consisting  of affine algebraic sets $A_n$ over $K$ and regular maps $F_{n m} \colon A_m \to A_n$ satisfying the following conditions:

\begin{enumerate}[(PSC-1)]
\item
$X_n$ is a constructible subset of $A_n$ for every $n \in \N$;
\item
$F_{n m}(X_m) \subset X_n$ and $f_{n m}$ is the restriction of $F_{n m}$ to $X_m$ for all $m,n \in \N$ such that  $m \geq n$.
\end{enumerate} 

The following  result is an essential ingredient in the proofs of Theorem \ref{t:closed-image} and Theorem \ref{t:reversible}.

 \begin{theorem}
 \label{t:proj-limit-constructible}
Let $K$ be an uncountable algebraically closed field and let   
  $(X_n,f_{n m})$ be a projective sequence of nonempty constructible sets over $K$.
   Then one has $ \varprojlim X_n \not= \varnothing$. 
\end{theorem}

Let us first prove Theorem \ref{t:proj-limit-constructible}
in the particular case where the projective sequence is given by inclusion maps 
(cf.   \cite[(CIP) p. 127]{gromov-esav}):

\begin{proposition}
\label{p:inter-constructible-not-empty}
Let $K$ be an uncountable algebraically closed field and let $A$ be an affine algebraic set 
over $K$.
Suppose that $(C_n)_{n \in \N}$ is a sequence of nonempty constructible subsets of $A$ such that 
$C_{n + 1} \subset C_n$ for all $n \in \N$.
Then one has $\bigcap_{n \in \N} C_n \not= \varnothing$.
\end{proposition}

We start by establishing two auxiliary results which are valid over any uncountable ground field.

\begin{lemma}
\label{l:intersection-countably-many-special-open}
Let $K$ be an uncountable (not necessarily algebraically closed) field and let $(Q_n)_{n \in \N}$ be a sequence of nonzero polynomials in $K[t_1,\dots,t_m]$.
Then there exists a point $a \in K^m$ such that $Q_n(a) \not= 0$ for all $n \in \N$.
\end{lemma}

\begin{proof}
We proceed by induction on $m$.
For $m = 1$, this follows from the fact that a nonzero polynomial  in one indeterminate has only finitely many zeroes and the fact that the union of a countable family of finite sets is countable.
Suppose now that $m \geq 2$ and that the result is true for polynomials in $m - 1$ indeterminates.
Let $S$ denote the set of $n \in \N$ such that the indeterminate $t_m$ occurs in $Q_n$.
Thus, we have $Q_n \in K[t_1,\dots,t_{m - 1}]$ for all $n \in \N \setminus S$.
For $n \in S$, let $R_n \in K[t_1,\dots,t_{m - 1}]$ denote the coefficient of the highest degree power of $t_m$ occurring in $Q_n$.
By our induction hypothesis, we can find   $b \in K^{m - 1}$ such that
$R_n(b) \not= 0$ for all $n \in S$ and $Q_n(b) \not= 0$ for all $n \in \N \setminus S$.
As $Q_n(b,t_m)$ is a nonzero polynomial in $t_m$ for all $n \in S$, it follows from the case  $m = 1$
that we can find $t \in K$    such that
$Q_n(b,t) \not= 0$ for all $n \in S$.
Then the point $a = (b,t) \in K^m$ satisfies $Q_n(a) \not=  0$ for all $n \in \N$.
 \end{proof}

\begin{lemma}
\label{l:intersection-countably-many-open}
Let $K$ be an uncountable (not necessarily algebraically closed) field and let 
$(\Omega_n)_{n \in \N}$ be a sequence of nonempty open subsets of  $K^m$.
Then one has $\bigcap_{n \in \N} \Omega_n \not= \varnothing$. 
\end{lemma}

\begin{proof}
As the special open subsets form a basis for the Zariski topology on $K^m$,
we can find, for each $n \in \N$, a nonzero polynomial $Q_n \in K[t_1,\dots,t_m]$ such that 
$V_n = K^m \setminus \Zer(Q_n) = \{a \in K^m : Q_n(a) \not= 0 \}$ satisfies  $V_n \subset \Omega_n$.
By Lemma \ref{l:intersection-countably-many-special-open}, we have $\bigcap_{n \in \N} V_n \not= \varnothing$.
Consequently, we also have $\bigcap_{n \in \N} \Omega_n \not= \varnothing$. 
\end{proof}

\begin{proof}[Proof of Proposition \ref{p:inter-constructible-not-empty}]
As the sequence of closed subsets $(\overline{C_n})_{n \in \N}$ is non-increasing, it eventually stabilizes. Thus, we can assume that $\overline{C_n} = A$ for all $n \in \N$.
\par
By Proposition \ref{p:constructible-contains-open}, we can find, for each $n \in \N$, a nonempty open subset $U_n$ of $A$ such that $U_n \subset C_n$. 
If $A = A_1 \cup A_2\cup \dots \cup A_s$ is the decomposition of $A$ into irreducible components, then
we have $U_n = (U_n \cap A_1) \cup (U_n \cap A_2) \cup \dots \cup (U_n \cap A_s) \not= \varnothing$.
It follows that we can find an index $1 \leq i \leq s$
and an increasing map $\varphi \colon \N \to \N$ such that
$U_{\varphi(n)} \cap A_i \not= \varnothing$ for all $n \in \N$.
  \par
   Let $d = \dim(A_i)$.
 Since $A_i$ is irreducible and the closed subset $F_n = A_i \setminus U_{\varphi(n)}$ is strictly contained in  $A_i$, we have 
$\dim(F_n) < d$.
On the other hand, it follows from Theorem \ref{t:finite-morphism} that  we can find a surjective finite morphism   
 $f \colon A_i \to K^d$.
 As every finite morphism is closed, the set $L_n = f(F_n)$ is closed in $K^d$.
We have  $\dim(L_n) \leq \dim(F_n) < d$ and therefore  $L_n \not= K^d$ for all $n \in \N$.
By Lemma \ref{l:intersection-countably-many-open}, the nonempty open subsets 
$\Omega_n = K^d \setminus L_n \subset K^d$ satisfy 
$\bigcap_{n \in \N} \Omega_n \not= \varnothing$.
As
$$
f(\bigcup_{n \in \N} F_n) = \bigcup_{n \in \N}f(F_n) = \bigcup_{n \in \N} L_n = K^d \setminus \bigcap_{n \in \N} \Omega_n ,
$$
it follows that
$$
f(\bigcup_{n \in \N} F_n) \not= K^d.
$$
As $f$ is surjective, this implies that $\bigcup_{n \in \N} F_n \not= A_i$ and hence 
$\bigcap_{n \in \N} U_{\varphi(n)} \not= \varnothing$.
Since   $U_{\varphi(n)} \subset C_{\varphi(n)} \subset C_n$ for all $n \in \N$, 
we conclude that  $\bigcap_{n \in \N} C_n \not= \varnothing$. 
\end{proof}

\begin{remark}
Proposition \ref{p:inter-constructible-not-empty} becomes  false  when the ground field $K$ is countable even if $K$ is algebraically closed (e.g., when $K$ is the algebraic closure of either $\Q$, or of the field $\F_p = \Z/p\Z$ of cardinality $p$ where $p$ is a prime number). Indeed, if $K$ is countable, say $K = \{a_n : n \in \N\}$, then the sequence of constructible subsets
$$
C_n = K \setminus \{a_0,a_1,\dots,a_n\} \subset K \quad (n \in \N)
$$
has an empty intersection.
\end{remark}

\begin{proof}[Proof of Theorem \ref{t:proj-limit-constructible}]
Let $(A_n,F_{n m})$ be a projective sequence of affine algebraic sets  and regular maps satisfying conditions (PSC-1) and (PSC-2) above. 
Let $(X_n',f_{n m}')$ denote the universal projective sequence associated with the projective sequence $(X_n,f_{n m})$.
For all $m \geq n$, the image set
  $f_{nm}(X_m) = F_{n m}(X_m)$ is a constructible subset of $A_n$ by Chevalley's theorem (Theorem \ref{t:chevalley}). 
 As the sequence  $f_{nm}(X_m)$, $m = n,n+1,\dots$, is a
non-increasing sequence of nonempty constructible subsets of the affine algebraic set $A_n$,
we deduce from Proposition \ref{p:inter-constructible-not-empty} that
$$
X_n' = \bigcap_{m \geq n}  f_{nm}(X_m) \not= \varnothing
$$
for all $n \in \N$.
 Thus, by Proposition \ref{p:universal}, it suffices to show that all maps $f_{n m}'$, $m \geq n$, are surjective.
 \par
 To see this, let $m,n \in \N$ with $m \geq n$ and
 suppose that  $x_n' \in X_n'$.
 Then, for all $k \geq n$, we have  $x_n' \in f_{n k}(X_k)$ so that we can find $y_k \in X_k$ such that $f_{n k}(y_k) = x_n'$.
 For $k \geq m$, the element $z_k = f_{m k}(y_k)$ satisfies 
 $f_{n m}(z_k) = f_{n m} \circ f_{m k}(y_k) = f_{n k}(y_k) = x_n'$.
 We deduce that $f_{n m}^{-1}(x_n') \cap f_{m k}(X_k) \not= \varnothing$ for all $k \geq m$.
 Now observe that $f_{n m}^{-1}(x_n') \cap f_{m k}(X_k)$ is constructible in $A_m$.
 Indeed,      
 $f_{n m}^{-1}(x_n') = F_{n m}^{-1}(x_n') \cap X_m$,  is constructible in $A_m$ since it is the intersection of a closed subset  with a constructible subset of $A_m$,  and  $f_{m k}(X_k) = F_{m k}(X_k)$   is constructible in $A_m$ by  Chevalley's theorem (Theorem \ref{t:chevalley}).
 By applying again Proposition \ref{p:inter-constructible-not-empty}, we deduce that
 $$
 f_{n m}^{-1}(x_n')  = \bigcap_{k \geq m} \left(f_{n m}^{-1}(x_n') \cap f_{m k}(X_k)\right) \not= \varnothing.
$$
Consequently,   the map $f_{n m}' \colon X_m' \to X_n'$ is surjective.
\end{proof}

\section{The closed image property}
\label{sec:closed-image}

This section contains the proofs of Theorem \ref{t:closed-image} 
and Corollary \ref{c:surjunctivity-rf}. 

\begin{proof}[Proof of Theorem \ref{t:closed-image}] 
Let $\tau \colon A^G \to A^G$ be an algebraic cellular automaton. Let $M \subset G$ be a memory set for $\tau$ and let $\mu \colon A^M \to A$ denote the associated local defining map.
\par
Suppose first that the group $G$ is countable.
Then we can find a sequence $(E_n)_{n \in \N}$ of finite subsets of $G$ such that $G = \bigcup_{n \in \N} E_n$,
$M \subset E_0$,  and $E_n \subset E_{n + 1}$ for all $n \in \N$.
Consider, for each $n \in \N$, the finite subset $F_n \subset G$ defined by  
$F_n = \{g \in G: gM \subset E_n\}$. 
Note that $G = \bigcup_{n \in \N} F_n$,
$1_G \in F_0$,  and $F_n \subset F_{n + 1}$
for all $n \in \N$.
\par
 It follows from \eqref{e;local-property} that if $x$ and $x'$ are elements in $A^G$ such that $x$ and $x'$ coincide on $E_n$  then the configurations $\tau(x)$ and $\tau(x')$ coincide on $F_n$.  Therefore, we can define a map
$\tau_n \colon A^{E_n} \to A^{F_n}$ by setting 
 $$
 \tau_n(u) = (\tau(x))\vert_{F_n}
 $$ 
 for all $u \in A^{E_n}$, where $x \in A^G$ denotes an arbitrary  configuration extending $u$.
Observe  that both $A^{E_n}$ and $A^{F_n}$ are affine algebraic sets as they are finite Cartesian powers of the affine algebraic set $A$. Moreover, it is clear from the fact that the map 
$\mu \colon A^M \to A$ is regular and formula \eqref{e;local-property}
that the map $\tau_n \colon A^{E_n} \to A^{F_n}$ is regular. 
\par
  Let now $y \in A^G$ and suppose that $y$ is in the closure of $\tau(A^G)$.
Then, for all $n \in \N$, we can find $z_n \in A^G$ such that 
\begin{equation}
\label{e;y-z-n-B-n}
 y\vert_{F_n} =   (\tau(z_n))\vert_{F_n}.
\end{equation}
Consider, for each $n \in \N$, the affine algebraic set $X_n \subset A^{E_n}$ defined by 
$X_n = \tau_n^{-1}( y\vert{F_n})$.
We have $X_n \not= \varnothing$ for all $n \in \N$ by \eqref{e;y-z-n-B-n}. 
Observe that, for all $m \geq n$, the restriction map $A^{E_m} \to A^{E_n}$ induces a regular  map
$f_{nm} \colon X_m \to X_n$. 
Conditions (PS-1) and (PS-2) are trivially satisfied so that
$(X_n,f_{nm})$ is a projective sequence of nonempty constructible (in fact, affine algebraic) sets.
By Theorem \ref{t:proj-limit-constructible}, we have $\varprojlim X_n \not= \varnothing$.
Choose an element $(x_n)_{n \in \N} \in \varprojlim X_n$. Thus
$x_n \in A^{E_n}$  and $x_{n + 1}$ coincides with $x_n$ on $E_n$  for all $n \in \N$. 
As $G = \cup_{n \in \N} E_n$, we deduce that there exists a (unique) configuration 
$x \in A^G$ such that $x\vert_{E_n} = x_n$ for all $n \in \N$. 
Moreover, we have $\tau(x)\vert_{F_n}= \tau_n(x_n) = y_n = y\vert_{F_n}$ for all $n$ since $x_n \in X_n$. 
As $G = \cup_{n \in \N} F_n$, this shows that $\tau(x) = y$.
This completes the proof  that $\tau$ has the closed image property in the case when $G$ is countable.
\par
Let us treat now the case of an arbitrary (possibly uncountable) group $G$.
Let $H$ denote the   subgroup of $G$ generated by $M$.
Observe that $H$ is countable since $M$ is finite.
The restriction cellular automaton $\tau_H \colon A^H \to A^H$ is algebraic by Theorem \ref{t:induction}.(iv).
Thus, by the first part of the proof, 
$\tau_H$ has the closed image property, that is,  
  $\tau_H(A^H)$ is closed in $A^H$ for the prodiscrete topology.  By applying Theorem 
  \ref{t:induction}.(iii), we deduce that $\tau(A^G)$ is also closed in $A^G$ for the prodiscrete topology.  Thus $\tau$ has the closed image property.
 \end{proof}

As the following example shows,
 Theorem \ref{t:closed-image} becomes false if the hypothesis saying that $K$ is algebraically closed is omitted.

\begin{example}
\label{ex:not-closed-image}
Take $K = \R$ and  consider the map  $\tau \colon \R^\Z \to \R^\Z$ defined by 
\begin{equation*}
\label{e:def-real-alg-ca-image-not-closed}
\tau(x)(n) = x(n + 1) - x(n)^2 \quad \text{ for all } x \in \R^\Z.
\end{equation*}
Clearly $\tau$ is an  algebraic cellular automaton over the group $\Z$ with memory set 
$M = \{0,1\}  $    and local defining map
$\mu \colon \R^M \to \R$ given by $\mu(x_0,x_1) = x_1 - x_0^2$ for all $(x_0,x_1) \in \R^2 = \R^M$.
 We claim that the image of  $\tau$  
is not closed in $\R^\Z$ for the prodiscrete topology.
 \par
Let us first show that $\tau(\R^{\Z})$ is dense in $\R^\Z$.
Let $y \in \R^\Z$ and let $F$ be a finite subset of $\Z$. 
Choose   $m \in \Z$ such that 
$F \subset [m,\infty)$.
Consider the configuration $x_F \in \R^\Z$ inductively defined by
$x_F(n) = 0$ for all $n \leq m$ and $x_F(n + 1) = y(n)  + x_F(n)^2$ for all $n \geq m$.  
We then have $\tau(x_F)(n) = y(n)$ for all $n  \geq m$  so that the configurations $\tau(x_F)$ and $y$ coincide on $[m,\infty)$ and hence on $F$.
Thus $y$ is in the closure of $\tau(\R^\Z)$.
\par
Consider now  the constant configuration $z \in \R^\Z$ defined by $z(n) = 1$  for all $n \in \Z$.
We claim that  the configuration $z$ is not in the image of $\tau$.
Suppose on the contrary that $z = \tau(x)$ for some $x \in \R^\Z$.
This means that $x(n + 1) = 1 + x(n)^2$ for all $n \in \Z$. It follows that $x(n) \geq 1$ and $x(n) < x(n + 1)$ for all $n \in \Z$
so that $x(n)$ must admit a finite limit as $n$ tends to $- \infty$. However, the existence of such a limit is impossible since the equation $\alpha = 1 + \alpha^2$ has no real roots.
This shows that $z$ is not in $\tau(\R^\Z)$. As $\tau(\R^\Z)$ is dense in $\R^\Z$, we conclude 
that $\tau(\R^\Z)$ is not closed in $\R^\Z$.
 \end{example}

\begin{remark}
\label{rem:not-closed-non-periodic}
More generally, if $G$ is any non-periodic group, 
then  one can construct an algebraic cellular automaton $\tau^G \colon \R^G  \to \R^G$ over the field $\R$ which does not have the closed image property.
Indeed, it suffices to choose an element of infinite order $g_0 \in G$ and consider the cellular automaton $\tau^G \colon \R^G \to \R^G$ obtained by induction from the cellular automaton 
$\tau \colon \R^\Z \to \R^\Z$ of the previous example, 
where we identify $\Z$ with  the subgroup of $G$ generated by $g_0$.
The fact that $\tau^G$ has the required properties follows from assertions (iii) and (iv) of Theorem \ref{t:induction}.
\end{remark}

Before proving Corollary \ref{c:surjunctivity-rf}, let us introduce additional notation.
 \par
Let   $A$, $M$, and $N$ be sets.
Suppose that we are given a map  $\rho \colon M \to N$.
Then $\rho$ induces a map $\rho^* \colon A^N \to A^M$ defined by $\rho^*(y) = y \circ \rho$ 
for all $y \in A^N$.

\begin{lemma}
\label{l:rho-induces-regular}
Let $K$ be a field and let $A$ be an affine algebraic set over $K$.
Suppose that we are given a map $\rho \colon M \to N$, where $M$ and $N$ are finite sets.
Then the induced map $\rho^* \colon A^N \to A^M$ is regular.
\end{lemma}

\begin{proof}
We have $\rho^*(y)(m) = y(\rho(m))$ for all $m \in M$ and $y \in A^N$. It follows that each coordinate map of $\rho^*$ is one of the projection maps $A^N \to A$ and is  therefore regular.
Consequently, $\rho^*$ is regular.
\end{proof}

Let $G$ be a group and let $A$ be a set.
Suppose that  $H$ is a subgroup of $G$.
Denote by  $\Fix(H)$ the subset of $A^G$ consisting of all   configurations
$x \in A^G$ which are fixed by $H$, that is, such that $hx = x$ for all $h \in H$.  
 Consider the set $H \backslash G = \{ Hg : g \in G \}$ consisting of all right cosets of $H$ in $G$ and the canonical surjection 
 $\rho_H \colon G \to H \backslash G$ which send each
$g \in G$ to $Hg$.
One immediately checks that $\rho_H^*(y) \in \Fix(H)$ for all $y \in A^{H \backslash G}$.
 Moreover, the map  $\rho_H^* \colon A^{H \backslash G} \to \Fix(H)$   is bijective 
 (see e.g. \cite[Proposition 1.3.3]{ca-and-groups-springer}).
Observe now that if $\tau \colon A^G \to A^G$ is a cellular automaton, then one has 
  $\tau(\Fix(H)) \subset \Fix(H)$ since $\tau$ is $G$-equivariant.
We  denote by   $\tau_H \colon \Fix(H) \to \Fix(H)$ the map  obtained by restriction of $\tau$,
and by $\widetilde{\tau_H} \colon A^{H \backslash G} \to A^{H \backslash G}$ the conjugate of 
$\tau$ by $\rho^*_H$, that is, the map given by
$\widetilde{\tau_H} = (\rho_H^*)^{-1} \circ \tau_H \circ \rho_H^*$.

\begin{proof}[Proof of Corollary \ref{c:surjunctivity-rf}]
Suppose that  $\tau \colon A^G \to A^G$ is an injective algebraic  cellular automaton.
Denote by $\FF$ the set of all finite index  subgroups of $G$.
\par
Let $H \in \FF$. Then $H \backslash G$ is finite.
 We claim that the map
$\widetilde{\tau_H} \colon A^{H \backslash G} \to A^{H \backslash G}$ is regular.
To see this, it suffices to prove that,
for each $ g \in G$, the map $ \pi_g \colon A^{H \backslash G} \to A$ defined by
$\pi_g(y) = \widetilde{\tau_H}(y)(Hg)$ is regular.
Choose a memory set  $M$  for $\tau$ and let $\mu \colon A^M \to A$ denote the associated local defining map.
Consider the map  $\psi \colon M \to H  \backslash G$  defined by $\psi(m) = \rho_H(gm)$ for all $m \in M$ and the induced map $\psi^* \colon A^{H \backslash G} \to A^M$.
Then we have $\pi_g = \mu \circ \psi^*$.
The map $\mu$ is regular since $\tau$ is algebraic. On the other hand, $\psi^*$ is regular by Lemma \ref{l:rho-induces-regular}.
It follows that $\pi_g$ is regular. this proves our claim.
  Now observe that $\tau_H \colon \Fix(H) \to \Fix(H)$ is injective since it is the restriction of $\tau$.  
As $\widetilde{\tau_H}$ is conjugate to $\tau_H$, we deduce that $\widetilde{\tau_H}$ is injective as well. It follows that $\widetilde{\tau_H}$ is  surjective by the Ax-Grothendieck theorem (Theorem \ref{t:ax-theorem}).
 Thus, $\tau_H$ is also surjective  and hence $ \Fix(H) = \tau_H(\Fix(H)) \subset \tau(A^G)$.
 \par
 Let $E \subset A^G$ denote the set of configurations whose orbit under the $G$-shift is finite. Then we have
$$
E = \bigcup_{H \in \FF} \Fix(H) \subset \tau(A^G).
$$
On the other hand, the residual finiteness of  $G$  implies that $E$ is dense in $A^G$ (see e.g. \cite[Theorem 2.7.1]{ca-and-groups-springer}).
As $\tau(A^G)$ is closed in $A^G$ by Theorem \ref{t:closed-image}, we conclude 
that $\tau(A^G) = A^G$.
 \end{proof}

\section{Reversibility}
\label{sec:reversibility}

\begin{proof}[Proof of Theorem \ref{t:reversible}]
 Let $\tau \colon A^G \to A^G$ be a bijective algebraic cellular automaton.
 We have to show that the inverse map $\tau^{-1} \colon A^G \to A^G$ is a cellular automaton. 
 \par 
Suppose first that the group $G$ is countable. 
 Let us show that the following local property is satisfied
by $\tau^{-1}$: 
\begin{itemize}
\item[($\ast$)]
there exists a finite subset $N \subset G$ such that, 
for any $y \in A^G$, the element $\tau^{-1}(y)(1_G)$ only depends on
the restriction of $y$ to $N$.
\end{itemize}
This will show that $\tau$ is reversible.
Indeed, if ($\ast$) holds for some finite subset $N \subset G$, then there exists a (unique) map 
$\nu \colon A^{N} \to A$ such that
$$
\tau^{-1}(y)(1_G) = \nu(y\vert_N)  
$$
for all $y \in A^G$.
Now, the $G$-equivariance of $\tau$ implies  the $G$-equivariance of its inverse map $\tau^{-1}$.
Consequently, we get
$$
\tau^{-1}(y)(g) = g^{-1}\tau^{-1}(y)(1_G) = \tau^{-1}(g^{-1}y)(1_G) =  \nu((g^{-1}y)\vert_N)
$$
for all $y \in A
^G$ and $g \in G$.
 which implies that $\tau^{-1}$ is the cellular automaton with memory set $N$ and local defining map $\nu$.
\par
Let us assume by contradiction that condition  ($\ast$) is not satisfied. 
Let $M$ be a memory set for $\tau$ such that $1_G \in M$.
Since $G$ is countable, we can find a sequence $(E_n)_{n \in \N}$ of finite subsets of $G$ such that $G = \bigcup_{n \in \N} E_n$,
$M \subset E_0$, and $E_n \subset E_{n + 1}$ for all $n \in \N$. 
 Consider, for each $n \in \N$, the finite subset $F_n \subset G$ defined by  
$F_n = \{g \in G: gM \subset E_n\}$. Note that $G = \bigcup_{n \in \N} F_n$,
$1_G \in F_0$,  and $F_n \subset F_{n + 1}$
for all $n \in \N$.
\par
Since    ($\ast$)  is not satisfied, we can find, for each $n \in \N$, two configurations 
$y_n', y_n'' \in A^G$ such that 
\begin{equation}
\label{e:property-y-n}
y_n'\vert_{F_n} = y_n''\vert_{F_n}
\quad \text{and} \quad 
\tau^{-1}(y_n')(1_G) \neq \tau^{-1}(y_n'')(1_G).
\end{equation}
Recall from the proof of Theorem \ref{t:closed-image}, 
that $\tau$ induces, for each $n \in \N$,  a regular map $\tau_n \colon A^{E_n} \to A^{F_n}$  given by $\tau_n(u) = (\tau(x))\vert_{F_n}$ for every $u \in A^{E_n}$, where $x \in A^G$ is any configuration extending $u$. 
  \par
Consider now, for each $n \in \N$, the subset $X_n \subset A^{E_n} \times A^{E_n}$ 
consisting of all pairs $(u,v)  \in A^{E_n} \times A^{E_n}$ such that 
$\tau_n(u) = \tau_n(v)$ and $u(1_G) \neq v(1_G)$. 
Note that $X_n$ is locally closed and hence constructible  in the affine algebraic set 
$A^{E_n} \times A^{E_n}$ for the Zariski topology since it is the intersection of a closed subset with an open subset.
Note also that $X_n$ is not empty since 
$$
((\tau^{-1}(y_n'))\vert_{E_n},(\tau^{-1}(y_n''))\vert_{E_n})   \in X_n
$$
 by \eqref{e:property-y-n}.
  Now observe that, for $m \geq n$, the restriction 
map  $\rho_{n m} \colon A^{E_m} \to A^{E_n}$
gives us a regular map 
$$
\pi_{n m}= \rho_{n m} \times \rho_{n m} \colon A^{E_n} \times A^{E_n} \to A^{E_n} \times A^{E_n}
$$
which induces by restriction a map $f_{n m} \colon X_m \to X_n$.
 Conditions (PS-1) and (PS-2) are trivially satisfied, so that
$(X_n,f_{nm})$ is a projective sequence of nonempty constructible sets.
 Thus,   we have $\varprojlim X_n \not= \varnothing$ by Theorem \ref{t:proj-limit-constructible}.
Choose an element $(p_n)_{n \in \N} \in \varprojlim X_n$. 
Thus
$p_n = (u_n,v_n) \in A^{E_n} \times A^{E_n}$  and $u_{n + 1}$ (resp. $v_{n + 1}$) coincides with $u_n$ (resp. $v_n$) on $E_n$  for all $n \in \N$. 
As $G = \cup_{n \in \N} E_n$, we deduce that there exists a (unique) configuration 
$x' \in A^G$ (resp. $x'' \in A^G$) such that $x'\vert_{E_n} = u_n$ 
(resp. $x''\vert_{E_n} = v_n$) for all $n \in \N$. 
Moreover, we have 
$$
(\tau(x'))\vert_{F_n}= \tau_n(u_n) = \tau_n(v_n) =  (\tau(x''))\vert_{F_n} 
$$
for all $n \in \N$.   
As $G = \cup_{n \in \N} F_n$, this shows that $\tau(x') = \tau(x'')$.
On the other hand, we have
$x'(1_G) = u_0(1_G) \not= v_0(1_G) = x''(1_G)$ and hence $x' \not= x''$.
   This contradicts the injectivity of $\tau$ and 
therefore completes the proof that $\tau$ is reversible in the case when $G$ is countable.
\par
 We now drop the countability assumption on $G$ and prove the theorem in its full  generality. Choose a memory set $M \subset G$ for $\tau$ and denote by $H$ the   subgroup of $G$ generated by  $M$.
Observe that $H$ is countable since $M$ is finite.
  By assertions (iv) and (i) of Theorem \ref{t:induction},
   the restriction cellular automaton $\tau_H \colon A^H \to A^H$ is algebraic and bijective. It then follows from the first part of the proof that $\tau_H$ is reversible.
 This implies that $\tau$ is reversible as well by assertion (ii) of Theorem \ref{t:induction}.
\end{proof}

\begin{remark}
Suppose that we are given a group $G$, a set $A$, and a bijective map $f \colon A \to A$.
Then the map $\tau \colon A^G \to A^G$ defined by $\tau(x)(g) = f(x(g))$ is a reversible cellular automaton with memory set $M = \{1_G\}$ and local defining map $f$.
The inverse cellular automaton $\tau^{-1} \colon A^G \to A^G$ is the cellular automaton with the same memory set $M$ and local defining map $f^{-1}$.
\par
By taking $A = K$, where  $K$ is an algebraically closed field of characteristic $p > 0$ and
 $f \colon K \to K$ the Frobenius automorphism
this gives an example of a reversible algebraic cellular automaton  $\tau \colon A^G \to A^G$ over $K$ whose inverse cellular automaton is not algebraic 
(cf. Remark~\ref{rem:ax}.(c)). 
\par
Similarly, by taking $A = \R$ and $f \colon \R \to \R$ given by $f(x) = x^3$, we get a reversible algebraic cellular automaton $\tau \colon A^G \to A^G$ over $\R$ whose inverse cellular automaton is not algebraic (cf. Remark \ref{rem:ax}.(d)).
\end{remark}
 
\begin{questions} 
The following questions are very natural:
\begin{enumerate}[(Q1)]
\item
Does there exist a bijective algebraic cellular automaton $\tau \colon A^\Z \to A^\Z$ over $\C$  whose inverse cellular automaton $\tau^{-1} \colon A^\Z \to A^\Z$ is not algebraic?
\item
Does there exist an injective algebraic cellular automaton $\tau \colon A^\Z \to A^\Z$ over $\R$ which is not surjective (cf. \cite[Remark (d) p.129]{gromov-esav})?
\end{enumerate}
\end{questions}

\bibliographystyle{siam}
\bibliography{algca}

\end{document}